\documentclass[11pt]{article}
\usepackage[utf8]{inputenc}
\usepackage[english]{babel}
\usepackage{amsthm}
\usepackage{graphicx}
\usepackage{url}
\usepackage{hyperref}
\usepackage{amsmath}
\usepackage{moreverb}
\usepackage{listings}
\usepackage{xcolor}

\begin{document}

\title{The recurrence formulas for primes and non-trivial zeros of the Riemann zeta function}

\author{Artur Kawalec}

\date{}
\maketitle

\begin{abstract}
In this article, we explore the Riemann zeta function with a perspective on primes and non-trivial zeros. We develop the Golomb's recurrence formula for the $n$th+1 prime, and assuming (RH), we propose an analytical recurrence formula for the $n$th+1 non-trivial zero of the Riemann zeta function. Thus all non-trivial zeros up the $n$th order must be known to generate the $n$th+1 non-trivial zero. We also explore a variation of the recurrence formulas for primes based on the prime zeta function, which will be a basis for the development of the recurrence formulas for non-trivial zeros based on the secondary zeta function. In the last part, we review the presented formulas and outline the duality between primes and non-trivial zeros. The proposed formula implies that all primes can be converted into an individual non-trivial zero (assuming RH), and conversely, all non-trivial zeros can be converted into an individual prime (not assuming RH). Also, throughout this article, we summarize numerical computation and verify the presented results to high precision.
\end{abstract}

\newpage
\section{Introduction}
The Riemann zeta function is defined by the infinite series
\begin{equation}\label{eq:20}
\zeta(s)=\sum_{n=1}^{\infty}\frac{1}{n^s},
\end{equation}
which is absolutely convergent for $\Re(s)>1$, where $s=\sigma+it$ is a complex variable. The values for the first few special cases are:

\begin{equation}\label{eq:9}
\begin{aligned}
\zeta(1) &\sim\sum_{n=1}^{k}\frac{1}{n}\sim\gamma+\log(k) \quad \text{as}\quad k\to \infty,\\
\zeta(2) &=\frac{\pi^2}{6}, \\
\zeta(3) &=1.20205690315959\dots, \\
\zeta(4) &=\frac{\pi^4}{90}, \\
\zeta(5) &=1.03692775514337\dots.
\end{aligned}
\end{equation}
For $s=1$, the series diverges asymptotically as $\gamma+\log(k)$, where $\gamma=0.5772156649\dots$ is the Euler-Mascheroni constant. The special values for even positive integer argument are given by the Euler's formula
\begin{equation}\label{eq:9}
\zeta(2k) = \frac{\mid B_{2k}\mid}{2(2k)!}(2\pi)^{2k},
\end{equation}
for which the value is expressed as a rational multiple of $\pi^{2k}$ where the constants $B_{2k}$ are Bernoulli numbers denoted such that $B_0=1$, $B_1=-\frac{1}{2}$, $B_2=\frac{1}{6}$ and so on. For odd positive integer argument, the values of $\zeta(s)$ converge to unique constants, which are not known to be expressed as a rational multiple of $\pi^{2k+1}$ as occurs in the even positive integer case. For $n=3$, the value is commonly known as  Ap\'ery's constant, who proved its irrationality.

At the heart of the Riemann zeta function are prime numbers, which are encoded by the Euler's product formula
\begin{equation}\label{eq:20}
\zeta(s)=\prod_{n=1}^{\infty}\left(1-\frac{1}{p_n^s}\right)^{-1}
\end{equation}
also valid for $\Re(s)>1$, where $p_1=2$, $p_2=3$, and $p_3=5$ and so on, denote the prime number sequence. The expression for the complex magnitude, or modulus, of the Euler prime product is

\begin{equation}\label{eq:7}
\mid \zeta(\sigma+it) \mid^2 = \frac{\zeta(4\sigma)}{\zeta(2\sigma)}\prod_{n=1}^\infty \left(1-\frac{\cos(t\log p_n)}{\cosh(\sigma\log p_n)}\right)^{-1}
\end{equation}
for $\sigma>1$, which for a positive integer argument $\sigma=k$ simplifies the zeta terms using (3), resulting in
\begin{equation}\label{eq:11}
\mid \zeta(k+it) \mid=(2\pi)^{k}\sqrt{\frac{|B_{4k}|(2k)!}{|B_{2k}|(4k)!}}\prod_{n=1}^\infty \left(1-\frac{\cos(t\log p_n)}{\cosh(k\log p_n)}\right)^{-1/2}.
\end{equation}
Using this form, the first few special values of this representation are

\begin{equation}\label{eq:29}
\begin{aligned}
\zeta(1) &\sim \frac{\pi}{\sqrt{15}}\prod_{n=1}^k \left(1-\frac{2}{p_n+p_n^{-1}}\right)^{-1/2} \sim e^{\gamma}\log(p_k),\\
\zeta(2) &= \frac{\pi^2}{\sqrt{105}}\prod_{n=1}^\infty \left(1-\frac{2}{p_n^2+p_n^{-2}}\right)^{-1/2}, \\
\zeta(3) &= \frac{\pi^3}{15} \sqrt{\frac{691}{3003}}\prod_{n=1}^\infty \left(1-\frac{2}{p_n^{3}+p_n^{-3}}\right)^{-1/2}, \\
\zeta(4) &= \frac{\pi^4}{45} \sqrt{\frac{3617}{17017}}\prod_{n=1}^\infty \left(1-\frac{2}{p_n^{4}+p_n^{-4}}\right)^{-1/2}, \\
\zeta(5) &= \frac{\pi^5}{225} \sqrt{\frac{174611}{323323}}\prod_{n=1}^\infty \left(1-\frac{2}{p_n^{5}+p_n^{-5}}\right)^{-1/2},
\end{aligned}
\end{equation}
where we let $t=0$ and reduced the hyperbolic cosine term as we have shown in [7][9]. The value for $\zeta(1)$ in terms of Euler prime product representation is asymptotic to $e^{\gamma}\log(p_k)$ due to Mertens's theorem as $k\to \infty$ [5][14]. Also, the arg of the Euler product can be found as
\begin{equation}\label{eq:20}
\text{arg } \zeta(\sigma+it) = -\sum_{n=1}^{\infty}\tan^{-1}\left(\frac{\sin(t\log p_n)}{p_n^{\sigma}-\cos(t\log p_n)}\right)
\end{equation}
thus writing the Euler product in polar form:
\begin{equation}\label{eq:20}
\zeta(s) = |\zeta(s)|e^{\text{i arg } \zeta(s)}.
\end{equation}
The Euler prime product permits the primes to be individually extracted from the infinite product under certain limiting conditions, as we have shown in [6], thus yielding the Golomb's formula for primes [4]. To illustrate this, when we expand the product we have

\begin{equation}\label{eq:20}
\zeta(s)=\left(1-\frac{1}{p_1^s}\right)^{-1}\left(1-\frac{1}{p_2^s}\right)^{-1}\left(1-\frac{1}{p_3^s}\right)^{-1}\ldots,
\end{equation}
and next, we wish to solve for the first prime $p_1$, then we have
\begin{equation}\label{eq:20}
p_1=\left(1-\frac{\epsilon_2(s)}{\zeta(s)}\right)^{-1/s},
\end{equation}
where
\begin{equation}\label{eq:20}
\epsilon_k(s)=\prod_{n=k}^{\infty}\left(1-\frac{1}{p_n^s}\right)^{-1}
\end{equation}
is the tail of Euler product starting at $p_k$. When we then consider the limit

\begin{equation}\label{eq:20}
p_1=\lim_{s\to \infty}\left(1-\frac{\epsilon_2(s)}{\zeta(s)}\right)^{-1/s},
\end{equation}
then $\epsilon_2(s)\to 1$ at a faster rate than the Riemann zeta function, that is $\zeta(s) \sim 1+O(p_1^{-s})$, while $\epsilon_2(s) \sim 1+O(p_2^{-s})$, and the gap $p_1^{-s}\gg p_2^{-s}$ is only widening as $s\to \infty$, hence the contribution due to Riemann zeta function dominates the limit, and the formula for the first prime becomes

\begin{equation}\label{eq:20}
p_1=\lim_{s\to \infty}\left(1-\frac{1}{\zeta(s)}\right)^{-1/s}.
\end{equation}
Numerical computation of (14) for $s=10$ and $s=100$ is summarized in Table $1$, where we observe convergence to $p_1$. The next prime in the sequence is found the same way by solving for $p_2$ in (10) to obtain
\begin{equation}\label{eq:20}
p_2=\lim_{s\to \infty}\left[1-\frac{\left(1-\frac{1}{p_1^s}\right)^{-1}\epsilon_3(s)}{\zeta(s)}\right]^{-1/s},
\end{equation}
where similarity as before, $\epsilon_3(s)\to 1$  at a faster rate than the Riemann zeta function and the contribution due to the first prime product $(1-p_1^{-s})^{-1}$ as $s\to\infty$, where it cancels the first prime product in $\zeta(s)$, so that $(1-p_1^{-s})\zeta(s) \sim 1+O(p_2^{-s})$, while $\epsilon_3(s) \sim 1+O(p_3^{-s})$, and the gap $p_2^{-s}\gg p_3^{-s}$ is increasing rapidly as $s\to \infty$, hence the contribution due to Riemann zeta function and the first prime product dominates the limit, and we have
\begin{equation}\label{eq:20}
p_2=\lim_{s\to \infty}\left[1-\frac{\left(1-\frac{1}{p_1^s}\right)^{-1}}{\zeta(s)}\right]^{-1/s}.
\end{equation}
Numerical computation of (16) for $s=10$ and $s=100$ is summarized in Table 1, and we observe convergence to $p_2$. And the next prime follows the same pattern $(1-p_1^{-s})(1-p_2^{-s})\zeta(s) \sim 1+O(p_3^{-s})$, while $\epsilon_4(s) \sim 1+O(p_4^{-s})$ which results in
\begin{equation}\label{eq:20}
p_3=\lim_{s\to \infty}\left[1-\frac{\left(1-\frac{1}{p_1^s}\right)^{-1}\left(1-\frac{1}{p_2^s}\right)^{-1}}{\zeta(s)}\right]^{-1/s}.
\end{equation}
Hence, this process continues for the $n$th+1 prime, and so if we define a partial Euler product up to the $n$th order as
\begin{equation}\label{eq:20}
Q_{n}(s)=\prod_{k=1}^{n}\left(1-\frac{1}{p_k^s}\right)^{-1}
\end{equation}
for $n>1$ and $Q_0(s)=1$, then we obtain the Golomb's formula for the $p_{n+1}$ prime
\begin{equation}\label{eq:20}
p_{n+1}=\lim_{s\to \infty}\left(1-\frac{Q_n(s)}{\zeta(s)}\right)^{-1/s}.
\end{equation}
We performed numerical computation of (19) in PARI/GP software package, as it is an excellent platform for performing arbitrary precision computations [10], and its functionality will be very useful for the rest of this article. Before running any script, we recommend to allocate alot of memory $\textbf{allocatemem(1000000000)}$, and setting precision to high value, for example $\textbf{\textbackslash p 2000}$. We tabulate the computational results in Table 1 for $s=10$ and $s=100$ case, and observe the convergence approaching to the $p_{n+1}$ prime based on the knowledge of all primes up to the $n$th order. When we compute for the $p_{1000}$ case, the $s=100$ variable is still too small to obverse correct convergence, hence we performed a very high precision computation for $n=999$ and $s=10000$ with precision set to $50000$ decimal places, and now the true value of the prime is revealed:

\begin{equation}\label{eq:20}
p_{1000}\approx7926.99958710978789301541492167\ldots.
\end{equation}
This formula will always converge because $p_n^{-s}\gg p_{n+1}^{-s}$ as $s\to \infty$, and also because the prime gaps are always bounded which will prevent higher order primes from modifying the main asymptote. It's just a matter of allowing the limit variable $s$ to tend a large value, however, as it seen it is not very practical for computing large primes, as very high arbitrary precision is required. The script in PARI is shown in Listing $1$ to compute the next prime using the Golomb's formula (19), which was used to compute Table $1$. The precision must be set very high, we generally set to $2000$ digits by default.
\begin{table}[hbt!]
\caption{The $p_{n+1}$ prime computed by equation (19) shown to $15$ decimal places.} 
\centering 
\begin{tabular}{c c c c} 
\hline\hline 
$n$ & $p_{n+1}$ & $s=10$ & $s=100$ \\[0.5ex] 
\hline 
$0$ & $p_1$ & 1.996546424130332  & 1.999999999999999 \\
$1$ & $p_2$ & 2.998128944738979 & 2.999999999999999 \\
$2$ & $p_3$ & 4.982816482987932 & 4.9999999999999991\\
$3$ & $p_4$ & 6.990872151877531 & 6.999999999999999 \\
$4$ & $p_5$ & 10.795904253794409 & 10.999999993885992 \\
$5$ & $p_6$ & 12.882858209904345 & 12.999999999999709 \\
$6$ & $p_7$ & 16.454690036492369 & 16.999997488242396 \\
$7$ & $p_8$ & 18.700432429563358 & 18.999999999042078 \\
$8$ & $p_9$ & 22.653649208924189 & 22.999999999980263 \\
$9$ & $p_{10}$ & 27.560268802131417 & 28.999632082761238  \\
$99$ & $p_{100}$ & 429.143320774398099 & 539.114941393037977 \\
$999$ & $p_{1000}$ & 5017.353999786395028 & 7747.370093956440561
\\ [1ex] 
\hline 
\end{tabular}
\label{table:nonlin} 
\end{table}

\newpage
\lstset{language=C,caption={PARI script for computing equation (19).},label=DescriptiveLabel,captionpos=b}
\begin{lstlisting}[frame=single]
\\ Define partial Euler product up to nth order
Qn(x,n)=
{
  prod(i=1,n,(1-1/prime(i)^x)^(-1));
}

\\ Compute the next prime
{
	n=10;	\\ set n
	s=100;  \\ set limit variable

        \\ compute next prime
	pnext=(1-Qn(s,n)/zeta(s))^(-1/s);
	print(pnext);
}
\end{lstlisting}

The Riemann zeta function has many representations. One common form is the alternating series representation

\begin{equation}\label{eq:20}
\zeta(s) = \frac{1}{1-2^{1-s}}\sum_{n=1}^{\infty} \frac{(-1)^{n+1}}{n^s},
\end{equation}
which is convergent for $\Re(s)>0$, with some exceptions at $\Re(s)=1$ due the constant factor. By the application of the Euler-Maclaurin summation formula, the main series (1) can also be extended to domain $\Re(s)>0$ by subtracting the pole in the limit as

\begin{equation}\label{eq:20}
\zeta(s)=\lim_{k\to \infty}\Big\{\sum_{n=1}^{k-1}\frac{1}{n^s}-\frac{k^{1-s}}{1-s}\Big\}.
\end{equation}
Equations (21) and (22) are hence valid in the critical strip region $0<\Re(s)<1$.

Another important representation of $\zeta(s)$ is the Laurent expansion about $s=1$ that gives a globally convergent series valid anywhere in the complex plane except at $s=1$ as

\begin{equation}\label{eq:20}
\zeta(s)=\frac{1}{s-1}+\sum_{n=0}^{\infty}\gamma_n\frac{(-1)^n(s-1)^n}{n!}.
\end{equation}
The coefficients $\gamma_n$ are the Stieltjes  constants, and $\gamma_0=\gamma$ is the usual Euler-Mascheroni constant. We observe that $\gamma_n$ are linear in the series, hence if we form a system of linear equations, then using the Cramer's rule and some properties of an Vandermonde matrix, we find that Stieltjes constants can be represented by a determinant of a certain matrix:

\begin{equation}\label{eq:20}
\gamma_n = \pm\det(A_{n+1})
\end{equation}
where the matrix $A_n(k)$ is matrix $A(k)$, but with an $n$th column swapped with a vector $B$ as given next

 \begin{gather}
A(k)= \begin{pmatrix}1 & -\frac{1}{1!} & \frac{1^2}{2!} & -\frac{1^3}{3!} &\dots & \frac{1^k}{k!}\\ 1 & -\frac{2}{1!} & \frac{2^2}{2!} & -\frac{2^3}{3!} & \dots & \frac{2^k}{k!} \\ 1 & -\frac{3}{1!} & \frac{3^2}{2!} & -\frac{3^3}{3!} & \dots & \frac{3^k}{k!}\\ \vdots & \vdots & \vdots & \vdots & \ddots  & \vdots\\ 1 & -\frac{(k+1)}{1!} & \frac{(k+1)^2}{2!} & -\frac{(k+1)^3}{3!} &\dots & \frac{(k+1)^k}{k!}\end{pmatrix}
\end{gather}
and

\begin{gather}
B(k)=
  \begin{pmatrix}
  \zeta(2)-1\\
  \zeta(3)-\frac{1}{2}\\
  \zeta(4)-\frac{1}{3}\\
  \vdots \\
  \ \zeta(k+1)-\frac{1}{k}\\
   \end{pmatrix}.
\end{gather}
The $\pm$ sign depends on $k$, but to ensure a positive sign, the size of $k$ must be a multiple of $4$.  Hence, the first few Stieltjes constants can be represented as:

\begin{gather}
\gamma_0= \lim_{k\to\infty}\begin{vmatrix} \zeta(2)-1 & -\frac{1}{1!} & \frac{1^2}{2!} & -\frac{1^3}{3!} &\dots & \frac{(-1)^k 1^k}{k!}\\ \zeta(3)-\frac{1}{2} & -\frac{2}{1!} & \frac{2^2}{2!} & -\frac{2^3}{3!} & \dots & \frac{(-1)^k 2^k}{k!} \\  \zeta(4)-\frac{1}{3} & -\frac{3}{1!} & \frac{3^2}{2!} & -\frac{3^3}{3!} & \dots & \frac{(-1)^k 3^k}{k!}\\ \vdots & \vdots & \vdots & \vdots & \ddots  & \vdots\\ \zeta(k+1)-\frac{1}{k} & -\frac{(k+1)}{1!} & \frac{(k+1)^2}{2!} & -\frac{(k+1)^3}{3!} &\dots & \frac{(-1)^k(k+1)^k}{k!}\end{vmatrix},
\end{gather}
and the next Stieltjes constant is
\begin{gather}
\gamma_1= \lim_{k\to\infty}\begin{vmatrix} 1 & \zeta(2)-1 & \frac{1^2}{2!} & -\frac{1^3}{3!} &\dots & \frac{(-1)^k 1^k}{k!}\\ 1 & \zeta(3)-\frac{1}{2} & \frac{2^2}{2!} & -\frac{2^3}{3!} & \dots & \frac{(-1)^k 2^k}{k!} \\  1 & \zeta(4)-\frac{1}{3} & \frac{3^2}{2!} & -\frac{3^3}{3!} & \dots & \frac{(-1)^k 3^k}{k!}\\ \vdots & \vdots & \vdots & \vdots & \ddots  & \vdots\\ 1 & \zeta(k+1)-\frac{1}{k}  & \frac{(k+1)^2}{2!} & -\frac{(k+1)^3}{3!} &\dots & \frac{(-1)^k(k+1)^k}{k!}\end{vmatrix},
\end{gather}
and the next is
\begin{gather}
\gamma_2=  \lim_{k\to\infty}\begin{vmatrix}1 & -\frac{1}{1!} &  \zeta(2)-1 & -\frac{1^3}{3!} &\dots & \frac{(-1)^k 1^k}{k!}\\ 1 & -\frac{2}{1!} & \zeta(3)-\frac{1}{2} & -\frac{2^3}{3!} & \dots & \frac{(-1)^k 2^k}{k!} \\  1 & -\frac{3}{1!} & \zeta(4)-\frac{1}{3}& -\frac{3^3}{3!} & \dots & \frac{(-1)^k 3^k}{k!}\\ \vdots & \vdots & \vdots & \vdots & \ddots  & \vdots\\ 1 & -\frac{(k+1)}{1!} & \zeta(k+1)-\frac{1}{k} & -\frac{(k+1)^3}{3!} &\dots & \frac{(-1)^k(k+1)^k}{k!}\end{vmatrix},
\end{gather}
and so on. In Table $2$, we compute the determinant formula (24) for the first $10$ Stieltjes constants for $k=500$, and observe the convergence.  In Listing $2$, the script in PARI to generate values for Table $2$ is also given. This shows that the $\gamma_n$ constants can be represented by $\zeta(n)$ at positive integer values as basis
\begin{equation}\label{eq:20}
\gamma_n = \lim_{k\to\infty}\Bigg\{C_{n,1}(k)+\sum_{m=2}^{k+1}(-1)^m C_{n,m}(k)\zeta(m)\Bigg\}
\end{equation}
where the expansion coefficients $C_{n,m}$ are rational and divergent, which grow very fast as $k$ increases.  The index $n\geq 0$ is the $n$th Stieltjes constant, and index $m\geq 1$ is for the $\zeta(m)$ basis value. These coefficients can be generated by expanding the determinant of $A_n$ using the Leibniz determinant rule along columns with the zeta values. For example, for $k=12$, which is a multiple of $4$, then the first few expansion coefficients are

\begin{equation}\label{eq:20}
\gamma_0 \approx -\frac{86021}{27720}+12\zeta(2)-66\zeta(3)+220\zeta(4)-495\zeta(5)+792\zeta(6)-\ldots
\end{equation}
The $C_{0,1}$ coefficient is the harmonic number $H_{12}$
\begin{equation}\label{eq:20}
C_{0,1}=-H_{k} = -\sum_{n=1}^{k}\frac{1}{n}
\end{equation}
and the next are
\begin{equation}\label{eq:20}
C_{0,m}=\binom{k}{m-1}.
\end{equation}
For the next $\gamma_n$, the first few coefficients for $k=12$ are

\begin{equation}\label{eq:20}
\gamma_1 \approx -\frac{1676701}{415800}+\frac{58301}{2310}\zeta(2)-\frac{72161}{420}\zeta(3)+\frac{76781}{126}\zeta(4)-\frac{79091}{56}\zeta(5)+\frac{80477}{35}\zeta(6)-\ldots,
\end{equation}
and for the next $\gamma_n$, we have

\begin{equation}\label{eq:20}
\gamma_2 \approx -\frac{5356117}{907200}+\frac{10418}{225}\zeta(2)-\frac{2270987}{6300}\zeta(3)+\frac{143644}{105}\zeta(4)-\frac{5520439}{1680}\zeta(5)+\frac{574108}{105}\zeta(6)-\ldots,
\end{equation}
and so on, but these coefficients are more difficult to determine and they diverge very fast.

\begin{table}[hbt!]
\caption{The first $30$ digits of $\gamma_n$ computed by equation (24) for $k=500$.} 
\centering 
\begin{tabular}{c c c} 
\hline\hline 
$n$ & $\gamma_n$ & Significant Digits \\ [0.5ex] 
\hline 
$0$  &  0.577215664901532860606512090082 & 34 \\
$1$  & -0.072815845483676724860586375874 & 34 \\
$2$  & -0.009690363192872318484530386035 & 33 \\
$3$  &  0.002053834420303345866160046542 & 32 \\
$4$  &  0.002325370065467300057468170177 & 31 \\
$5$  &  0.00079332381730106270175333487\underline{7} & 30 \\
$6$  & -0.0002387693454301996098724218\underline{4}2 & 29 \\
$7$  & -0.0005272895670577510460740975\underline{0}7 & 29 \\
$8$  & -0.000352123353803039509602052\underline{1}77 & 28 \\
$9$  & -0.000034394774418088048177914\underline{6}91 & 28 \\
$10$  & 0.0002053328149090647946837\underline{2}1922 & 26
\\ [1ex] 
\hline 
\end{tabular}
\label{table:nonlin} 
\end{table}

\newpage
\lstset{language=C,caption={PARI script for computing equation (24)},label=DescriptiveLabel,captionpos=b}
\begin{lstlisting}[frame=single]
{
    n = 0;    \\ set nth Stieltjes constant
    k = 100;  \\ set limit variable

    An=matrix(k,k); \\ allocate matrix

    \\ load matrix An
    for(j=1,k,
        for(i=1,k,
            if(j==1+n,An[i,j]=zeta(i+1)-1/i,
               An[i,j]=(-i)^(j-1)/factorial(j-1))));

    \\ compute determinant of An
    yn = matdet(An);
    print(yn);
}
\end{lstlisting}

The Hadamard infinite product formula is another global analytically continued representation of (1) to the whole complex plane

\begin{equation}\label{eq:20}
\zeta(s)=\frac{\pi^{\frac{s}{2}}}{2(s-1)\Gamma(1+\frac{s}{2})}\prod_{\rho}^{}\left(1-\frac{s}{\rho}\right)
\end{equation}
having a simple pole at $s=1$, and at the heart of this form is an infinity of complex non-trivial zeros $\rho_n=\sigma_n+it_n$, which are constrained to lie in the critical strip $0<\Re(s)<1$ region. The infinite product is assumed to be taken over zeros in conjugate pairs. Hardy proved that there is an infinity of non-trivial zeros on the critical line at $\sigma=\frac{1}{2}$. It is not yet known whether there are non-trivial zeros off of the critical line in the range $0<\Re(s)<1$ other than $\sigma=\frac{1}{2}$, a problem of the Riemann Hypothesis (RH). To date, there has been a very large number of zeros verified numerically to lie on the critical line, and none was ever found off of the critical line. The first few non-trivial zeros on the critical line $\rho_n=\frac{1}{2}+it_n$ have imaginary components  $t_1 = 14.13472514...$, $t_2 = 21.02203964...$, $t_3 = 25.01085758...$ which were originally found numerically using a solver, but if (RH) is true, then can be computed analytically by the formula presented later in this article. Also, we will interchangeably refer to $\rho_n$ or $t_n$ to imply a non-trivial zero.

The Hadamard product representation can be interpreted as a volume of an s-ball (that is for a ball of complex dimension $s$). For a positive integer $n$, the n-ball defines all points satisfying $\Omega=\{x_1^2+x_2^2+x_3^2\dots +x_n^2\leq R^n\}$, and integrating gives the total volume

\begin{equation}\label{eq:20}
V(n)=\underset{\Omega}{\int\int\int\ldots\int} dx_1 dx_2 dx_3\ldots dx_n=K(n)R^n,
\end{equation}
where
\begin{equation}\label{eq:20}
K(n) = \frac{\pi^{\frac{n}{2}}}{\Gamma(1+\frac{n}{2})}
\end{equation}
is the proportionality constant. Now, when generalizing the n-ball to an s-ball of complex $s$ dimension for $\zeta(s)$, we can identify that the terms involving $\pi$ and $\Gamma(s)$ function is $K(s)$, and that the radius of the s-ball is the remaining product involving the non-trivial zeros
\begin{equation}\label{eq:25}
R(s)^s = \frac{1}{2(s-1)}\prod_{\rho}^{}\left(1-\frac{s}{\rho}\right)
\end{equation}
which is actually the Riemann xi function $\xi(s)$ multiplied by $(s-1)^{-1}$. Thus
\begin{equation}\label{eq:25}
\zeta(s) = V_s = K(s)R(s)^s
\end{equation}
can be understood as a volume quantity, which when packed into an s-ball, then the radius function in this form is being described by explicitly the non-trivial zeros. The trivial zeros at negative even integers $-2,-4,-6\ldots -2n$ are then the zeros of the proportionality constant due to the pole of $\Gamma(s)$. For example, if we consider $s=2$, then

\begin{equation}\label{eq:25}
\begin{aligned}
\zeta(2) & = K(2)R(2)^2 \\
         & =  \pi R^2
\end{aligned}
\end{equation}
where $R=\sqrt{\pi/6}=0.7236012545\ldots$ is the radius to give the volume quantity for $\zeta(2)$, which from (1) can be understood as packing the areas of squares with $1/n$ sides into a circle. And similarly for $s=3$
\begin{equation}\label{eq:25}
\begin{aligned}
\zeta(3) & =  K(3)R(3)^3 \\
         & =  \frac{4}{3}\pi R^3
\end{aligned}
\end{equation}
where $R=0.6595972037\ldots$ is the radius to give the volume quantity for Ap\'ery's constant $\zeta(3)$, which from (1) can be understood as packing the volumes of cubes with $1/n$ sides into a sphere. Hence in this view, the non-trivial zeros are governing the radius quantity of an s-ball, essentially encoding the volume information of $\zeta(s)$, and while the trivial zeros are just the zeros of the proportionality constant $K(s)$, which has a role of scaling the values of non-trivial zeros across the dimension $s$ to the values that they currently are, and perhaps even on the critical line. If we plot the radius in the range $1<\sigma<\infty$, we find a local minima for $R$ which occurs between $s=2$ and $s=3$ at $s_{min}= 2.8992592006...$ and $R_{min}=0.6592484066\ldots$. That would mean that the s-ball would reach minimum radius $R_{min}$ at $s_{min}$.

Furthermore, if we consider the complex magnitude for $\zeta(s)$ for representations (21) and (22), and note that at each non-trivial zero on the critical line, a harmonic series is induced from which we can obtain formulas for the Euler-Mascheroni constant $\gamma$ expressed as a function of a single non-trivial on the critical line zero as
\begin{equation}\label{eq:20}
\gamma = \lim_{k\to \infty}\Bigg\{2\sum_{v=1}^{k}\sum_{u=v+1}^{k}\frac{(-1)^{u}(-1)^{v+1}}{\sqrt{uv}}\cos(t^{}_n \log(u/v))-\log(k)\Bigg\}
\end{equation}
and the second formula as
\begin{equation}\label{eq:20}
\gamma = \lim_{k\to \infty}\Bigg\{\frac{k+1}{(\frac{1}{2})^2+t^{2}_{n}}-2\sum_{v=1}^{k}\sum_{u=v+1}^{k}\frac{1}{\sqrt{uv}}\cos(t^{}_n \log(u/v))-\log(k)\Bigg\},
\end{equation}
where it is assumed the index variables satisfy $u>v$ starting with $v=1$ as we have shown in [8][9]. Thus, any individual non-trivial zero on the critical line $t_n$ can be converted to $\gamma$, which is independent on (RH). As a numerical example, for $t_1$ and $k=10^5$, we compute $\gamma=0.5772\underline{1}81648\ldots$ accurate to $5$ decimal places, however, the computation becomes more difficult as it grows as $O(k^2)$ due to the double series. And if we subtract equations (43) and (44), then we obtain a relation

\begin{equation}\label{eq:20}
\frac{1}{|\rho_n|^2}=\frac{1}{(\frac{1}{2})^2+t^{2}_{n}} = \lim_{k\to\infty}\frac{2}{\sqrt{k}}\sum_{m=1}^{k}\frac{1}{\sqrt{m}}\cos(t_n\log(m/k))
\end{equation}
whereby any individual non-trivial zero can be converted to its absolute value on the critical line.  Also next, the infinite sum over non-trivial zeros

\begin{equation}\label{eq:20}
\sum_{n=1}^{\infty}\frac{1}{|\rho_n|^2}=\frac{1}{2}\gamma+1-\frac{1}{2}\log(4\pi),
\end{equation}
is an example of secondary zeta function family which will be discussed later.

There is also another whole side to the theory of the Riemann zeta function concerning the prime counting function $\pi(n)$ up to a given quantity $n$, and the non-trivial zero counting function $N(T)$ up to a given quantity T. It is natural to take the logarithm of the Euler prime product yielding a sum

\begin{equation}\label{eq:20}
\log[\zeta(s)] = \sum_{n=1}^{\infty}\sum_{m=1}^{\infty}\frac{1}{m}\frac{1}{p_n^{ms}}=s\int_{0}^{\infty}J(x)x^{-s-1}dx \quad \Re(s)>1
\end{equation}
from which motivates to define a step function $J(x)$ that increases by $1$ at each prime, by $\frac{1}{2}$ at prime square, by $\frac{1}{3}$ at prime cubes, and so on, as shown in [3, p.22] and [15]. Riemann then expressed $J(x)$ by Fourier inversion as

\begin{equation}\label{eq:20}
J(x)=\frac{1}{2\pi i}\int_{a-i\infty}^{a+i\infty}\log[\zeta(s)]\frac{x^s}{s}ds \quad(a>1).
\end{equation}
After finding a suitable expansion for $\log[\zeta(s)]$ in terms of non-trivial zeros using the xi function Weierstrass product over non-trivial zeros
\begin{equation}\label{eq:25}
\xi(s) = \frac{1}{2}\prod_{\rho}^{}\left(1-\frac{s}{\rho}\right),
\end{equation}
with its relation to the zeta by $\xi(s)=\pi^{-\frac{s}{2}}\zeta(s)(s-1)\Gamma(1+\frac{s}{2})$, then after a very detailed and a lengthy analysis as shown in Edwards [3], the main formula for $J(x)$ appears as

\begin{equation}\label{eq:20}
J(x) = \text{Li}(x)-\sum_{\rho}^{}\text{Li}(x^{\rho})-\log(2)+\int_{x}^{\infty}\frac{dt}{t(t^2-1)\log(t)}
\end{equation}
for $x>1$, where Li$(x)$ is a logarithmic integral, and then by applying M\"{o}bius inversion leads to recovering
\begin{equation}\label{eq:20}
\pi(x)=\sum_{n=1}^{\infty}\frac{\mu(n)}{n}J(x^{1/n}).
\end{equation}
Hence, through this formula, the non-trivial zeros are shown to be involved in the generation of primes. Although applying M\"{o}bius inversion in (51) to recover $\pi(n)$ is somewhat circular in this case, because one needs to have knowledge of all primes by $\mu(n)$, however, the main prime content is still in $J(x)$, which comes from the contribution of non-trivial zero terms. In Fig 1, we plot $J(x)$ using (50) and observe the curve approach the step function as more non-trivial zeros (taken in conjugate-pairs) are used. 

Furthermore, in analysis by LeClair [11] concerning $N(T)$, it is found that $n$th non-trivial zeros satisfy the following transcendental equation:

\begin{equation}\label{eq:20}
\frac{t_n}{2\pi}\log\left(\frac{t_n}{2\pi e}\right)+\lim_{\delta\to 0}\frac{1}{\pi}\text{arg } \zeta(\frac{1}{2}+it_n+\delta)=n-\frac{11}{8},
\end{equation}
however, the contribution to due to arg function is very small, and only provides fine level tuning to the overall equation, hence when dropping the arg term, LeClair obtained an approximate asymptotic formula for non-trivial zeros via the Lambert function $W(x)e^{W(x)}=x$ transformation:

\begin{equation}\label{eq:20}
t_{n}\approx 2\pi\frac{n-\frac{11}{8}}{W\left(\frac{n-\frac{11}{8}}{e}\right)}.
\end{equation}
It turns out that this approximation works very well with an accuracy down to a decimal place. For example, with this formula, we can quickly approximate a $10^{100}$ zero:

\begin{equation}\label{eq:20}
\begin{aligned}
t_{10^{100}}\approx &&  28069038384289406990319544583825640008454803016284\\
&& 6045192360059224930922349073043060335653109252473.23351
\end{aligned}
\end{equation}
in less than one second, and it  should be accurate to within a decimal place. The Lambert function can be computed efficiently for large input argument, and the approximated values for $t_n$ get better for higher zeros as $n\to\infty$. In fact, LeClair computed the largest non-trivial zero known to date for $n=10^{10^6}$ using this method [12].

Also, very little is known about the properties of non-trivial zeros. For example, they are strongly believed to be simple, but remains unproven. And in the works by Wolf [16], a large sample of non-trivial zeros was numerically expanded into continued fractions, from which it was possible to compute the Khinchin’s constant, which strongly suggests they are irrational.

In this article, we propose an analytical recurrence formula for $t_{n+1}$, very similar to the Golomb's formula for primes, thus all non-trivial zeros up to $t_n$ must be known in order to compute the $t_{n+1}$ zero. The formula is based on a certain representation of the secondary zeta function

\begin{equation}\label{eq:20}
Z(s) = \sum_{n=1}^{\infty}\frac{1}{t_{n}^{s}}
\end{equation}
in the works of Voros [13], for $s>1$, which is not involving non-trivial zeros, thus avoiding circular reasoning. There is alot of work already on the secondary zeta functions published in the literature, especially concerning the meromorphic extension of $Z(s)$ via the Mellin transform techniques and tools of spectral theory.

We now introduce the main result of this paper. Assuming (RH), the full recurrence formula for the $t_{n+1}$ non-trivial zero is:

\begin{equation}\label{eq:20}
t_{n+1} = \lim_{m\to\infty}\left[\frac{(-1)^{m+1}}{2}\left(\frac{1}{(2m-1)!}\log (|\zeta|)^{(2m)}\big(\frac{1}{2}\big)+\sum_{k=1}^{\infty}\frac{1}{\left(\frac{1}{2}+2k\right)^{2m}}-2^{2m}\right)-\sum_{k=1}^{n}\frac{1}{t_{k}^{2m}}\right]^{-\frac{1}{2m}}
\end{equation}
for $n\geq 0$, thus all non-trivial zeros up the $n$th order must be known in order to generate the $n$th+1 non-trivial zero. This formula is a solution to

\begin{equation}\label{eq:20}
\zeta(s)=0
\end{equation}
where $s=\rho_n=\frac{1}{2}+it_n$ for $\sigma_n=\frac{1}{2}$, and the zeros $t_n$ are real and ordered $0<t_1<t_2<t_3<\ldots t_{n}$. This formula is satisfied by all representations of $\zeta(s)$ on the critical strip, such as by (21), (22), (23), (36), and so on. And in the next sections, we will develop this formula, and explore some its variations, and then we will numerically compute non-trivial zeros to high precision. We will also discuss some possible limitations to this formula for $n\to \infty$.

In the last section, we will discuss formulas for $t_{n}$ which actually can be related to the primes themselves, and that one could compute $t_n$ as a function of all primes. And conversely, one could compute any individual prime $p_{n}$ as a function of all non-trivial zeros.

\section{A variation of the $n$th+1 prime formula}
Golomb described several variations of the prime formulas of the form (19), one such is

\begin{equation}\label{eq:20}
p_{n+1}=\lim_{s\to \infty}\left[\zeta(s)-Q_{n}(s)\right]^{-1/s},
\end{equation}
which will serve to motivate the next result, which is based on the prime zeta function, and that will then serve as a basis for the development of an analogue formula for the $n$th+1 non-trivial zero formula in the next section.

The prime zeta function is an analogue of (1), but instead of summing over reciprocal integer powers, we sum over reciprocal prime powers as
\begin{equation}\label{eq:20}
P(s) = \sum_{n=1}^{\infty}\frac{1}{p_{n}^{s}}.
\end{equation}
When we consider the expanded sum
\begin{equation}\label{eq:20}
P(s) = \frac{1}{p_{1}^{s}}+\frac{1}{p_{2}^{s}}+\frac{1}{p_{3}^{s}}+\ldots
\end{equation}
then similarly as before, we wish to solve for $p_1$, and obtain

\begin{equation}\label{eq:20}
\frac{1}{p_{1}^{s}}=P(s) -\frac{1}{p_{2}^{s}}-\frac{1}{p_{3}^{s}}-\ldots
\end{equation}
which leads to
\begin{equation}\label{eq:20}
p_1=\left(P(s) -\frac{1}{p_{2}^{s}}-\frac{1}{p_{3}^{s}}-\ldots\right)^{-1/s}.
\end{equation}
If we then consider the limit,
\begin{equation}\label{eq:20}
p_1=\lim_{s\to\infty}\left(P(s) -\frac{1}{p_{2}^{s}}-\frac{1}{p_{3}^{s}}-\ldots\right)^{-1/s}
\end{equation}
then we find that the higher order primes decay faster than $P(s)$, namely, $P(s)\sim p_1^{-s}$, while the tailing error is $O(p_2^{-s})$, and so $P(s)$ dominates the limit. Since $p_1^{-s}\gg p_2^{-s}$, hence we have
\begin{equation}\label{eq:20}
p_1=\lim_{s\to\infty}\left[P(s)\right]^{-1/s}.
\end{equation}
To find $p_2$ we consider (60) again
\begin{equation}\label{eq:20}
p_2=\lim_{s\to\infty}\left[P(s)-\frac{1}{p_1^s}-\frac{1}{p_3^s}\ldots\right]^{-1/s},
\end{equation}
and when taking the limit, then we must keep $p_1$, while the higher order primes decay faster, namely, $P(s)-p_1^{-s}\sim p_2^{-s}$, while the tailing error is $O(p_3^{-s})$, and so $P(s)-p_1^{-s}$ dominates the limit. Since $p_2^{-s}\gg p_3^{-s}$, hence we have
\begin{equation}\label{eq:20}
p_2=\lim_{s\to\infty}\left[P(s)-\frac{1}{p_1^s}\right]^{-1/s}.
\end{equation}
And similarly, the next prime is found the same way, but this time we must retain the two previous primes
\begin{equation}\label{eq:20}
p_3=\lim_{s\to\infty}\left[P(s)-\frac{1}{p_1^s}-\frac{1}{p_2^s}\right]^{-1/s}.
\end{equation}
Hence in general, if we define a partial prime zeta function up to the $n$th order
\begin{equation}\label{eq:20}
P_n(s) = \sum_{k=1}^{n}\frac{1}{p_{k}^{s}},
\end{equation}
then the $n$th+1 prime is
\begin{equation}\label{eq:20}
p_{n+1}=\lim_{s\to\infty}\left[P(s)-P_{n}(s)\right]^{-1/s}.
\end{equation}
At this point, knowing $P(s)$ by the original definition (60) leads to circular reasoning, hence we seek to find other representations for $P(s)$ that don't involve primes directly. We explore the well-known relation

\begin{equation}\label{eq:20}
\log[\zeta(s)]=\sum_{k=1}^{\infty}\frac{P(ks)}{k}
\end{equation}
and then by applying M\"{o}bius inversion leads to

\begin{equation}\label{eq:20}
P(s)=\sum_{k=1}^{\infty}\mu(k)\frac{\log(ks)}{k},
\end{equation}
where $\mu(k)$ is the M\"{o}bius function, which however, still depends on the primes just like (51) for $J(x)$, so it may not be best candidate for $P(s)$. And if there are other representations for $P(s)$ not involving primes, then one could certainly use them, but we are unaware of such. But to verify (69), we pre-compute $P(s)$ using primes to high precision instead, thus introducing circular reasoning, Hence, we pre-compute $P(s)$ for $s=10$ and $s=100$ as
\begin{equation}\label{eq:20}
P(10)=9.936035744369802178558507001 \times 10^{-4}\ldots
\end{equation}
and
\begin{equation}\label{eq:20}
P(100)=7.888609052210118073520537827\times 10^{-31}\ldots
\end{equation}
using a neat remainder estimation technique of (71) developed by Cohen in [2]. Next, we summarize computation for $p_{n+1}$ by the recurrence formula (69) in Table 3, and observe the convergence to the $p_{n+1}$ prime, just as the Golomb's formula for primes. And as before, the convergence works because

\begin{equation}\label{eq:20}
O(p_n^{-s})\gg O(p_{n+1}^{-s})\quad \text{as }s\to\infty,
\end{equation}
and also that the prime gaps are bounded, which prevents any higher order primes from modifying the main asymptote.

\begin{table}[hbt!]
\caption{The $p_{n+1}$ prime computed by equation (69) shown to $15$ decimal places.} 
\centering 
\begin{tabular}{c c c c} 
\hline\hline 
$n$ & $p_{n+1}$ & $s=10$ & $s=100$ \\[0.5ex] 
\hline 
$0$ & $p_1$ & 1.996543079767713  & 1.999999999999999 \\
$1$ & $p_2$ & 2.998128913153986 & 2.999999999999999 \\
$2$ & $p_3$ & 4.982816481260483 & 4.999999999999999\\
$3$ & $p_4$ & 6.990872151845387 & 6.999999999999999 \\
$4$ & $p_5$ & 10.79590425378718 & 10.999999993885992 \\
$5$ & $p_6$ & 12.88285820990352 & 12.999999999999709 \\
$6$ & $p_7$ & 16.45469003649213 & 16.999997488242396 \\
$7$ & $p_8$ & 18.70043242956331 & 18.999999999042078 \\
$8$ & $p_9$ & 22.65364920892418 & 22.999999999980263 \\
$9$ & $p_{10}$ & 27.5602688021314 & 28.999632082761238
\\ [1ex] 
\hline 
\end{tabular}
\label{table:nonlin} 
\end{table}

\newpage
\section{The recurrence formula for non-trivial zeros}
The secondary zeta function has been studied in the literature, and there has been interesting developments concerning the analytical extension to the whole complex plane for

\begin{equation}\label{eq:20}
Z(s) = \sum_{n=1}^{\infty}\frac{1}{t_{n}^{s}}
\end{equation}
which has many parallels with the zeta function. In this article, the symbol $Z$ is implied, and is not related to the Hardy-Z function. For the first few special values the $Z(s)$ yields

\begin{equation}\label{eq:9}
\begin{aligned}
Z(2) &=\frac{1}{2}(\log |\zeta|)^{(2)}\big(\frac{1}{2}\big)+\frac{1}{8}\pi^2+\beta(2)-4 \\
     &=0.023104993115418\dots, \\
     &\\
Z(3) &= 0.00072954\dots, \\
     &\\
Z(4) &=-\frac{1}{12}(\log |\zeta|)^{(4)}\big(\frac{1}{2}\big)-\frac{1}{24}\pi^4-4\beta(4)+16 \\
     &= 0.00037172599285\dots, \\
     &\\
Z(5) &= 0.00000223\dots.
\end{aligned}
\end{equation}

The special values for even positive integer argument $Z(2m)$ is:
\begin{equation}\label{eq:20}
\begin{aligned}
Z(2m) = (-1)^m \bigg[-\frac{1}{2(2m-1)!}(\log |\zeta|)^{(2m)}\big(\frac{1}{2}\big)+\\
         -\frac{1}{4}\left[(2^{2m}-1)\zeta(2m)+2^{2m}\beta(2m)\right]+2^{2m}\bigg]
\end{aligned}
\end{equation}
and is found in [13,p. 693] by works of Voros, and it's originally denoted as $\mathcal{Z}(2\sigma)$. This formula is a sort of an analogue for Euler's formula (3) for $\zeta(2n)$, and is valid for $m\geq 1$, where $m$ is an integer, and $\beta(s)$ is the Dirichlet beta function

\begin{equation}\label{eq:20}
\beta(s) = \sum_{n=0}^{\infty}\frac{(-1)^n}{(2n+1)^s}=\prod_{n=1}^{\infty}\left(1-\frac{\chi_{4}(p_n)}{p_n^s}\right)^{-1},
\end{equation}
where $\chi_4$ is the Dirichlet character modulo $4$. The value for $\beta(2)$ is the Catalan's constant. In (76), the odd values for $Z(2m+1)$ were computed numerically by summing $25000$ zeros, as it is not known whether there is a closed-form representation similarly as for the $\zeta(2m+1)$ case, and so the given values could only be accurate to several decimal places.  The formula (77) assumes (RH), and is a result of a complicated development to meromophically extend (75) to the whole complex plane using tools from spectral theory.
Furthermore, using the relation, also found in [13,p. 681] as
\begin{equation}\label{eq:20}
\frac{1}{2^s}\zeta\big(s,\frac{5}{4}\big)=\sum_{k=1}^{\infty}\frac{1}{\left(\frac{1}{2}+2k\right)^s}=2^s\left[\frac{1}{2}\left((1-2^{-s})\zeta(s)+\beta(s)\right)-1\right],
\end{equation}
from which we have several variations of (77) for $Z(2m)$ as
\begin{equation}\label{eq:20}
Z(2m) = \frac{(-1)^{m+1}}{2} \left[\frac{1}{(2m-1)!}\log (|\zeta|)^{(2m)}\big(\frac{1}{2}\big)+\sum_{k=1}^{\infty}\frac{1}{\left(\frac{1}{2}+2k\right)^{2m}}-2^{2m}\right]
\end{equation}
and another as
\begin{equation}\label{eq:20}
Z(2m) = \frac{(-1)^{m+1}}{2} \left[\frac{1}{(2m-1)!}\log (|\zeta|)^{(2m)}\big(\frac{1}{2}\big)+\frac{1}{2^{2m}}\zeta\big(2m,\frac{5}{4}\big)-2^{2m}\right].
\end{equation}
The expressions involving the $\log (|\zeta|)^{(2m)}\big(\frac{1}{2}\big)$ term can be computed numerically and independently of the non-trivial zeros, and there is no known closed-form representation of it, but there is for the odd values
\begin{equation}\label{eq:20}
\log (|\zeta|)^{(2m+1)}\big(\frac{1}{2}\big)=\frac{1}{2}(2m)!(2^{2m+1}-1)\zeta(2m+1)+\frac{1}{4}\pi^{2m+1}|E_{2m}|,
\end{equation}
where $E_{2m}$ are Euler numbers [13,p. 686]. Unfortunately, the $\log (|\zeta|)^{(2m+1)}(\frac{1}{2})$ term is not involved in the computation of $Z(m)$ for $m>1$. Also, the infinite series in (80) is related to the Hurwitz zeta function, and it can also be separated into two parts involving the zeta function and the beta function, which can then be related to primes via the Euler product, which we will come back to shortly.

Now we will follow the same program that we did for the prime zeta function as outlined in equations (58) to (69). If we begin with the secondary zeta function
\begin{equation}\label{eq:20}
Z(s) = \frac{1}{t_{1}^{s}}+\frac{1}{t_{2}^{s}}+\frac{1}{t_{3}^{s}}+\ldots
\end{equation}
and then solving for $t_1$  we obtain

\begin{equation}\label{eq:20}
\frac{1}{t_{1}^{s}}=Z(s) -\frac{1}{t_{2}^{s}}-\frac{1}{t_{3}^{s}}-\ldots
\end{equation}
and then we get
\begin{equation}\label{eq:20}
t_1=\left(Z(s) -\frac{1}{t_{2}^{s}}-\frac{1}{t_{3}^{s}}-\ldots\right)^{-1/s}.
\end{equation}
If we then consider the limit
\begin{equation}\label{eq:20}
t_1=\lim_{s\to\infty}\left(Z(s) -\frac{1}{t_{2}^{s}}-\frac{1}{t_{3}^{s}}-\ldots\right)^{-1/s}
\end{equation}
then, since $O(Z(s))\sim O(t_1^{-s})$, and so the higher order non-trivial zeros decay as $O(t_2^{-s})$ faster than $Z(s)$, and so $Z(s)$ dominates the limit, hence we have
\begin{equation}\label{eq:20}
t_1=\lim_{s\to\infty}\left[Z(s)\right]^{-1/s}.
\end{equation}
Now, substituting representation (80) for $Z(s)$ into (87), and $s$ is now assumed be an integer as a limit variable $2m$, then we get a direct formula for $t_1$ as

\begin{equation}\label{eq:20}
t_{1} = \lim_{m\to\infty}\left[\frac{(-1)^{m+1}}{2}\left(\frac{1}{(2m-1)!}\log (|\zeta|)^{(2m)}\big(\frac{1}{2}\big)+\sum_{k=1}^{\infty}\frac{1}{\left(\frac{1}{2}+2k\right)^{2m}}-2^{2m}\right)\right]^{-\frac{1}{2m}}.
\end{equation}
Next we numerically verify this formula in PARI, and the script is shown in Listing $3$. We broke up the representation (88) into several parts A to D. Also, sufficient memory must be allocated and precision set to high before running the script. We utilize the Hurwitz zeta function representation, since it is available in PARI, and the \textbf{derivnum} function for computing the $m$th derivative very accurately for high $m$. The results are summarized in Table $4$ for various limit values of $m$ from low to high, and we can observe the convergence to the real value as $m$ increases. Already at $m=10$ we get several digits of $t_1$, and at $m=100$ we get over $30$ digits. We performed even higher precision computations, and the result is clearly converging to $t_1$.

\begin{table}[hbt!]
\caption{The computation of $t_1$ by equation (88) for different $m$.} 
\centering 
\begin{tabular}{c c c} 
\hline\hline 
m & $t_1$ (First 30 Digits)  & Significant Digits\\ [0.5ex] 
\hline 
$1$ & 6.578805783608427637281793074245 & 0  \\
$2$ & 12.806907343833847091925940068962 & 0 \\
$3$ & 13.809741306055624728153992726341 & 0 \\
$4$ & 14.038096225961619450676758199577 & 0 \\
$5$ & 14.\underline{1}02624784431488524304946186056 & 1 \\
$6$ & 14.\underline{1}23297656314161936112154413740 & 1 \\
$7$ & 14.1\underline{3}0464459254236820197453483721 & 2 \\
$8$ & 14.1\underline{3}3083993992268169646789606564 & 2 \\
$9$ & 14.13\underline{4}077755601528384660110026302 & 3 \\
$10$ & 14.13\underline{4}465134057435907124435534843 & 3 \\
$15$ & 14.1347\underline{2}1950874675119831881762569 & 5 \\
$20$ & 14.13472\underline{5}096741738055664458081219 & 6\\
$25$ & 14.13472514\underline{1}055464326339414131271 & 9 \\
$50$ & 14.134725141734693\underline{7}89641535771021 & 16 \\
$100$ & 14.134725141734693790457251983562 & 34
\\ [1ex] 
\hline 
\end{tabular}
\label{table:nonlin} 
\end{table}

\newpage

\lstset{language=C,caption={PARI script for computing equation (88).},label=DescriptiveLabel,captionpos=b}
\begin{lstlisting}[frame=single]
{
    \\ set limit variable
    m = 250;

    \\ compute parameters A to D
    A = derivnum(x=1/2,log(zeta(x)),2*m);	
    B = 1/factorial(2*m-1);
    C = 2^(2*m);
    D = (2^(-2*m))*zetahurwitz(2*m,5/4);

    \\ compute Z(2m)
    Z = (-1)^(m+1)*(1/2)*(A*B-C+D);

    \\ compute t1
    t1 = Z^(-1/(2*m));
    print(t1);
}
\end{lstlisting}

Next, we perform a higher precision computation for $m=250$ case, and the result is

\begin{equation}\label{eq:20}
\begin{aligned}
t_1=14.13472514173469379045725198356247027078425711569924 & \\
     317568556746014996342980925676494901\underline{0}212214333747\ldots
\end{aligned}
\end{equation}
accurate to $87$ decimal places. In order to find the second non-trivial zero, we comeback to (83), and solving for $t_2$ yields

\begin{equation}\label{eq:20}
t_2=\lim_{s\to\infty}\left(Z(s) -\frac{1}{t_{1}^{s}}-\frac{1}{t_{3}^{s}}-\ldots\right)^{-1/s}
\end{equation}
and since the higher order zeros decay faster than $Z(s)-t_1^{-s}$, we then have

\begin{equation}\label{eq:20}
t_2=\lim_{s\to\infty}\left(Z(s) -\frac{1}{t_{1}^{s}}\right)^{-1/s}
\end{equation}
and the zero becomes

\begin{equation}\label{eq:20}
t_{2} = \lim_{m\to\infty}\left[\frac{(-1)^{m+1}}{2}\left(\frac{1}{(2m-1)!}\log (|\zeta|)^{(2m)}\big(\frac{1}{2}\big)+\sum_{k=1}^{\infty}\frac{1}{\left(\frac{1}{2}+2k\right)^{2m}}-2^{2m}\right)-\frac{1}{t_{1}^{2m}}\right]^{-\frac{1}{2m}}.
\end{equation}
A numerical computation for $m=250$ yields

\begin{equation}\label{eq:20}
\begin{aligned}
t_2=21.0220396387715549926284795938969027773\underline{3}355195796311 & \\
     4759442381621433519190301896683837161904986197676\ldots
\end{aligned}
\end{equation}
which is accurate to $38$ decimal places, and we assumed $t_1$ used was already pre-computed to $2000$ decimal places by other means. We cannot use the same $t_1$ computed earlier with same precision, as it will cause a self-cancelation in (91), and so, the numerical accuracy of $t_{n}$ must be much higher than $t_{n+1}$ to guarantee convergence. And continuing on, the next zero is computed as
\begin{equation}\label{eq:20}
t_{3} = \lim_{m\to\infty}\left[\frac{(-1)^{m+1}}{2}\left(\frac{1}{(2m-1)!}\log (|\zeta|)^{(2m)}\big(\frac{1}{2}\big)+\sum_{k=1}^{\infty}\frac{1}{\left(\frac{1}{2}+2k\right)^{2m}}-2^{2m}\right)-\frac{1}{t_{1}^{2m}}-\frac{1}{t_{2}^{2m}}\right]^{-\frac{1}{2m}}.
\end{equation}
A numerical computation for $m=250$ yields

\begin{equation}\label{eq:20}
\begin{aligned}
t_3=25.010857580145688763213790992562821818659549\underline{6}5846378 & \\
     3317371101068278652101601382278277606946676481041\ldots
\end{aligned}
\end{equation}
which is accurate to $43$ decimal places, and we assumed $t_1$ and $t_2$ was used to high enough precision which was $2000$ decimal places in this example.
Hence, just like for the $n$th+1 Golomb prime recurrence formulas and the prime zeta function $P(s)$, the same limit works for non-trivial zeros. As a result, if we define a partial secondary zeta function up to the $n$th order
\begin{equation}\label{eq:20}
Z_n(s) = \sum_{k=1}^{n}\frac{1}{t_{k}^{s}},
\end{equation}
then the $n$th+1 non-trivial zero is

\begin{equation}\label{eq:20}
t_{n+1}=\lim_{m\to\infty}\left[Z(m)-Z_{n}(m)\right]^{-1/m}
\end{equation}
and the main recurrence formula:

\begin{equation}\label{eq:20}
t_{n+1} = \lim_{m\to\infty}\left[\frac{(-1)^{m+1}}{2}\left(\frac{1}{(2m-1)!}\log (|\zeta|)^{(2m)}\big(\frac{1}{2}\big)+\sum_{k=1}^{\infty}\frac{1}{\left(\frac{1}{2}+2k\right)^{2m}}-2^{2m}\right)-\sum_{k=1}^{n}\frac{1}{t_{k}^{2m}}\right]^{-\frac{1}{2m}}.
\end{equation}
One can actually use any number of representations for $Z(s)$, and the challenge will be find more efficient algorithms to compute them. And finally, we report a numerical result for $Z(2m)$ for $m=250$ as:
\begin{equation}\label{eq:20}
\begin{aligned}
Z = 7.18316934899718140841650578011166023417090863769600 & \\
     8517536818521464413577481501771580460474425539208\times 10^{-576}\ldots.
\end{aligned}
\end{equation}
From this number, we extracted the first $10$ non-trivial zeros, which are summarized in Table $5$ for $k=250$. The previous non-trivial zeros used were already known to high precision to $2000$ decimal places in order to compute the $t_{n+1}$ zero. One cannot use the same $t_n$ obtained earlier because it will cause self-cancelation in (98), and the accuracy for $t_n$ must be much higher than $t_{n+1}$ to ensure convergence. Initially we started with an accuracy of $87$ digits after decimal place for $t_1$, and then it dropped to $7$ to $12$ digits by the time it gets to $t_{10}$ zero. There is also a sudden drop in accuracy when the gaps get too small. Hence, these formulas are not very practical for computing high zeros as large numerical precision is required, for example, at the first Lehmer pair at $t_{6709}=7005.06288$, the gap between next zero is about $\sim 0.04$.  Also, the average gap between zeros gets smaller as $t_{n+1}-t_{n}\sim\frac{2\pi}{\log(n)}$, making the use of this formula progressively harder and harder to compute.

\begin{table}[hbt!]
\caption{The $t_{n+1}$ computed by equation (98).} 
\centering 
\begin{tabular}{c c c c} 
\hline\hline 
$n$ & $t_{n+1}$ & $m=250$ & Significant Digits \\ [0.5ex] 
\hline 
$0$ & $t_{1}$ & 14.134725141734693790457251983562 & 87 \\
$1$ & $t_{2}$  & 21.022039638771554992628479593896 & 38 \\
$2$ & $t_{3}$ & 25.010857580145688763213790992562 & 43 \\
$3$ & $t_{4}$  & 30.424876125859513\underline{2}09940851142395 & 16 \\
$4$ & $t_{5}$  & 32.9350615877391896906623689640\underline{7}3 & 29 \\
$5$ & $t_{6}$  & 37.58617815882567125\underline{7}190902153280 & 18 \\
$6$ & $t_{7}$  & 40.918719012147\underline{4}63977678179889317 & 13 \\
$7$ & $t_{8}$  & 43.3270732809149995194961\underline{1}7449701 & 22 \\
$8$ & $t_{9}$  & 48.005150\underline{8}79831498066163921378664 & 7 \\
$9$ & $t_{10}$ & 49.77383247767\underline{2}299146155484901550 & 12
\\ [1ex] 
\hline 
\end{tabular}
\label{table:nonlin} 
\end{table}

\section{Duality between primes and non-trivial zeros}
We outline the duality between primes and non-trivial zeros. The Golomb's recurrence formula (19) is an exact formula for the $n$th+1 prime

\begin{equation}\label{eq:20}
p_{n+1}=\lim_{s\to \infty}\left(1-\frac{Q_n(s)}{\zeta(s)}\right)^{-1/s},
\end{equation}
and the Hadamard product formula establishes $\zeta(s)$ as a function of non-trivial zeros:

\begin{equation}\label{eq:20}
\zeta(s)=\frac{\pi^{\frac{s}{2}}}{2(s-1)\Gamma(1+\frac{s}{2})}\prod_{\rho}^{}\left(1-\frac{s}{\rho}\right).
\end{equation}
Hence, this is a pathway from non-trivial zeros to the primes and without assuming (RH), as the Hadamard product is over all zeros. On the other hand, the recurrence formula for the $n$th+1 non-trivial zero is
\begin{equation}\label{eq:20}
\begin{aligned}
t_{n+1} =  \lim_{m\to\infty}&\Bigg[\frac{(-1)^{m+1}}{2}\Big(\frac{1}{(2m-1)!}\log (|\zeta|)^{(2m)}\big(\frac{1}{2}\big)-2^{2m+1}+\\
          & +2^{2m-1}\big((1-2^{-2m})\zeta(2m)+\beta(2m)\big)\Big)-\sum_{k=1}^{n}\frac{1}{t_k}\Bigg]^{-\frac{1}{2m}}
\end{aligned}
\end{equation}
where now one could substitute the Euler product for the zeta and beta functions, or both, which is what we will do next. We have

\begin{equation}\label{eq:20}
\left(1-2^{-2m}\right)\zeta(2m)=\prod_{n=2}^{\infty}\left(1-\frac{1}{p_n^{2m}}\right)^{-1}
\end{equation}
and
\begin{equation}\label{eq:20}
\beta(2m)=\prod_{n=2}^{\infty}\left(1-\frac{\chi_4(p_n)}{p_n^{2m}}\right)^{-1}.
\end{equation}
As a result,

\begin{equation}\label{eq:20}
\begin{aligned}
t_{n+1} =  \lim_{m\to\infty}&\Bigg[\frac{(-1)^{m+1}}{2}\Bigg(\frac{1}{(2m-1)!}\log (|\zeta|)^{(2m)}\big(\frac{1}{2}\big)-2^{2m+1}+\\
          & +2^{2m-1}\Big(\prod_{n=2}^{\infty}\big(1-p_n^{-2m}\big)^{-1}+\prod_{n=2}^{\infty}\big(1-\chi_4(p_n)p_n^{-2m}\big)^{-1}\Big)\Bigg)-\sum_{k=1}^{n}\frac{1}{t_k}\Bigg]^{-\frac{1}{2m}}
\end{aligned}
\end{equation}
which completes the pathway from primes to non-trivial zeros. We note that these formulas are independent, and thus avoid any circularity, however, the recurrence formula for $t_{n+1}$ is dependent on (RH). And finally, in Appendix A, we present a PARI script to compute (105) recursively for several zeros.

\section{Conclusion}
We explored various representations of the Riemann zeta function, such as the Euler prime product, the Laurent expansion, and the Golomb's recurrence formula for primes. The Golomb's formula is a basis for developing similar recurrence formulas for the $n$th+1 non-trivial zeros via an independent formula for the secondary zeta function $Z(2m)$, which does not involve non-trivial zeros. Hence, the non-trivial zeros can be extracted under the right excitation in the limit, just like prime numbers. We verified these formulas numerically, and they indeed do converge to $t_{n+1}$. The difficultly lies in computation of the $\log (|\zeta|)^{(2m)}(\frac{1}{2})$ term. We utilized the PARI/GP software package for computing $Z(2m)$ for $m=250$, and the first zero $t_1$ achieves $87$ correct digits after the decimal place. Presently, computing beyond that caused the test computer to run out of memory. And so, if better and more efficient methods for computing $Z(2m)$ are developed, then more higher zeros can be computed accurately. But even then, computing up to a millionth zero for example, would be almost insurmountable. The only open question is whether the recurrence for the non-trivial zeros will hold up, namely the limit $O(t_n^{-s})\gg O(t_{n+1}^{-s})$ as $s\to\infty$, as the average gap between non-trivial zeroes decreases $t_{n+1}-t_{n}\sim\frac{2\pi}{\log(n)}$ as $n\to\infty$. In case of the Golomb's formula for primes, this gap is bounded.

These formulas also suggest a new criterion for (RH). It suffices to take a first zero $t_1$ represented by (88) which depends on (RH) as
\begin{equation}\label{eq:20}
t_{1} = \lim_{m\to\infty}\left[\frac{(-1)^{m+1}}{2}\left(\frac{1}{(2m-1)!}\log (|\zeta|)^{(2m)}\big(\frac{1}{2}\big)+\sum_{k=1}^{\infty}\frac{1}{\left(\frac{1}{2}+2k\right)^{2m}}-2^{2m}\right)\right]^{-\frac{1}{2m}}
\end{equation}
and passing it through to any number of representations of $\zeta(s)$ valid in the critical strip to work out
\begin{equation}\label{eq:20}
\zeta(\frac{1}{2}+i t_1)=0.
\end{equation}
For example, if we take equation (45) and substitute $t_1$ as
\begin{equation}\label{eq:20}
\frac{1}{|\rho_1|^2}=\frac{1}{(\frac{1}{2})^2+t_1^{2}} = \lim_{k\to\infty}\frac{2}{\sqrt{k}}\sum_{m=1}^{k}\frac{1}{\sqrt{m}}\cos(t_1\log(m/k)),
\end{equation}
then recovering $t_1$ would imply (RH) if there was a way work out the series.
We also would like to extend these formulas for the secondary beta function

\begin{equation}\label{eq:20}
B(s) = \sum_{n=1}^{\infty}\frac{1}{r_{n}^{s}},
\end{equation}
where $r_n$ are imaginary components of non-trivial zeros of $\beta(s)$ on the critical line. For example, the first few zeros are $r^{}_1 = 6.02094890...$, $r^{}_2 = 10.24377030...$, $r^{}_3 = 12.98809801...$. Then, the proposed recurrence formula would be
\begin{equation}\label{eq:20}
r_{n+1}=\lim_{s\to\infty}\left[B(s)-B_{n}(s)\right]^{-1/s},
\end{equation}
where
\begin{equation}\label{eq:20}
B_n(s) = \sum_{n=1}^{n}\frac{1}{r_{n}^{s}}
\end{equation}
is the partial secondary beta function up to the $n$th order. And just like for the Dirichlet beta, the same could potentially apply to other Dirichlet L-functions.

Finally, we highlighted the duality between primes and non-trivial zeros where it is possible convert non-trivial zeros into an individual prime, and conversely, to convert all primes into an individual non-trivial zero.

\texttt{Email: art.kawalec@gmail.com}

\newpage
\section{Appendix A}
The script in Listing $4$ computes the $n$th+1 non-trivial recursively from a set of primes by equation (105). The parameter $\textbf{pmax}$ specifies the number of primes to use for the Euler product. The starting limiting variable is $\textbf{m}$, and at each iteration $\textbf{m}$ is decreased by a pre-set amount $\textbf{step\textunderscore m}$, so that the accuracy for $t_n$ will be greater than for $t_{n+1}$ in order to avoid self-cancelation. The values for computed zeros are stored in an array, and the partial secondary zeta $Z_n$ is computed at every iteration. By leveraging these parameters, the output can converge to different values, and in some cases will not converge. We optimized them to give $4$ zeros accurately, and beyond that it doesn't converge and then $m$ has to be increased to a larger value. The results of running this script are summarized in Table $6$. As before, we obtain $t_1$ accurate to $87$ decimal places, but $t_2$ now is accurate to $26$ decimal places, and the next zero to $12$ and $1$ decimal places respectively. At this point the iteration has ran its course. We would like to increase $m$, but presently is outside the range of the test computer.

\begin{table}[hbt!]
\caption{The $t_{n+1}$ by PARI script in Listing $4$} 
\centering 
\begin{tabular}{c c c c c} 
\hline\hline 
$m$ & $n$ & $t_{n+1}$ & First $30$ digits of computed results & Significant Digits \\ [0.5ex] 
\hline 
250 & $0$  & $t_{1}$ & 14.134725141734693790457251983562 & 87 \\
175 & $1$  & $t_{2}$ & 21.0220396387715549926284795\underline{9}4245 & 26 \\
100 & $2$  & $t_{3}$ & 25.01085758014\underline{5}177681574221575793 & 12 \\
25 & $3$  & $t_{4}$  & 30.\underline{4}13415903597141481192661667214 & 1
\\ [1ex] 
\hline 
\end{tabular}
\label{table:nonlin} 
\end{table}
\newpage

\lstset{language=C,caption={PARI script for generating non-trivial zeros from primes.},label=DescriptiveLabel,captionpos=b}
\begin{lstlisting}[frame=single][hbt!]
{
    m = 250;          \\ starting limit variable m
    step_m = -75;     \\ decrease limit step_m	
    pmax = 2000;      \\ set max number of primes
    tn = vector(100); \\ allocate vector to hold zeros
    n=1;              \\ init non-trivial zero counter

    \\ start loop
    while(m != 0,

        \\ compute parameters A to D
        A = derivnum(x=1/2,log(zeta(x)),2*m);	
        B = 1/(factorial(2*m-1));
        C = 2^(2*m+1); D = 2^(2*m-1);

        \\ compute Euler products
        P1 = prod(i=2,pmax,(1-1/prime(i)^(2*m))^(-1));
        P2 = prod(i=2,pmax,
        (1-(-1)^((prime(i)-1)/2)/prime(i)^(2*m))^(-1));

        \\ compute Z(2m)
        Z = 0.5*(-1)^(m+1)*(A*B-C+D*(P1+P2));

        \\ compute Zn up to nth-1 order
        if(n==1,Zn=0,
            for(j=1,n-1,Zn = Zn + 1/tn[j]^(2*m)));

        \\ compute and print tn
        tn[n] = (Z-Zn)^(-1/(2*m));
        print(m, ":", tn[n]);

        m = m+step_m;   \\ decrease m by step_m
        n = n+1;        \\ increment zero counter
    )
}
\end{lstlisting}

\end{document}